\input ./style/arxiv-vmsta.cfg
\documentclass[numbers,compress,v1.0.1]{vmsta}
\usepackage{vtexbibtags}

\volume{5}% Updated by VTEXPTS2LaTeX.exe, 18.12.2018 09:32
\issue{4}% Updated by VTEXPTS2LaTeX.exe, 18.12.2018 09:32
\pubyear{2018}
\firstpage{521}% Updated by VTEXPTS2LaTeX.exe, 18.12.2018 09:32
\lastpage{541}% Updated by VTEXPTS2LaTeX.exe, 18.12.2018 09:32
\aid{VMSTA106}
\doi{10.15559/18-VMSTA106}% Updated by VTEXPTS2LaTeX.exe, 14.06.2018
\articletype{research-article}

\usepackage{authorquery}
\usepackage{amsfonts}
\usepackage{amsmath}
\usepackage{amsthm,amssymb}
\numberwithin{equation}{section}
\allowdisplaybreaks
\usepackage{longtable}
\theoremstyle{definition}

\startlocaldefs
\newcommand{\di}{\mathrm{d}}
\newcommand{\balpha}{\boldsymbol\alpha}

\newcommand{\f}{\frac}
\newcommand{\omin}{{\varOmega}_i^n}
\newcommand{\p}{\partial}
\newcommand{\pt}{{\mathcal P}_t}
\newcommand{\px}{{\mathcal P}_x}
\newcommand{\gm}{{\mathcal G}}
\newcommand{{\pip}}{p_{i+1}^n}
\newcommand{{\ppi}}{p_i^n}
\newcommand{{\pim}}{p_{i-1}^n}
\newcommand{{\xim}}{x_{i-1/2}}
\newcommand{{\xip}}{x_{i+1/2}}
\newcommand{{\xpm}}{x_{i \pm1/2}}

\newcommand{\lleft}{\left}
\newcommand{\rright}{\right}
\allowdisplaybreaks
\endlocaldefs

\begin{document}

\begin{frontmatter}
\pretitle{Research Article}

\title{Stable L\'evy diffusion and related model fitting}

\author[a]{\inits{P.}\fnms{Paramita}~\snm{Chakraborty}\thanksref{cor1}\ead[label=e1]{chakrabp@stat.sc.edu}\orcid{0000-0003-0021-6955}}
\thankstext[type=corresp,id=cor1]{Corresponding author.}
\author[b]{\inits{X.}\fnms{Xu}~\snm{Guo}\ead[label=e2]{guoxu1014@hkbu.edu.hk}}
\author[c]{\inits{H.}\fnms{Hong}~\snm{Wang}\ead[label=e3]{hwang@math.sc.edu}}

\address[a]{Department of Statistics, \institution{University of South Carolina}, Columbia, SC, \cny{USA}}
\address[b]{Department of Mathematics, \institution{Hong Kong Baptist University}, \cny{Hong Kong}}
\address[c]{Department of Mathematics, \institution{University of South Carolina}, Columbia, SC, \cny{USA}}

\markboth{P. Chakraborty et al.}{Stable L\'evy diffusion and related model fitting}

%\author[]{\inits{}\fnms{}~\snm{}\thanksref{f1}\ead[label=e1]{}}
%\author[]{\inits{}\fnms{}~\snm{}\thanksref{f1}\thanksref{cor1}
%\ead[label=e2]{}}
%\thankstext[type=corresp,id=cor1]{Corresponding author.}
%\address[]{\institution{}, ..., \cny{}}
%\address[]{\institution{}, ..., \cny{}}

%\thankstext[id=f1]{}

%\dedicated{}

%\markboth{Authors}{Title}
%\markboth{}{}

\begin{abstract}
A fractional advection-dispersion equation (fADE) has been advocated for
heavy-tailed flows where the usual Brownian diffusion models fail. A
stochastic differential equation (SDE) driven by a stable L\'evy process
gives a forward equation that matches the space-fractional
advection-dispersion equation and thus gives the stochastic framework of
particle tracking for heavy-tailed flows. For constant advection and
dispersion coefficient functions, the solution to such SDE itself is a
stable process and can be derived easily by least square parameter fitting
from the observed flow concentration data. However, in a more generalized
scenario, a closed form for the solution to a stable SDE may not exist. We
propose a numerical method for solving/generating a stable SDE in a
general set-up. The method incorporates a discretized finite volume scheme
with the characteristic line to solve the fADE or the forward equation for
the Markov process that solves the stable SDE. Then we use a numerical
scheme to generate the solution to the governing SDE using the fADE
solution. Also, often the functional form of the advection or dispersion
coefficients are not known for a given plume concentration data to start
with. We use a Levenberg--Marquardt (L-M) regularization method to estimate
advection and dispersion coefficient function from the observed data (we
present the case for a linear advection) and proceed with the SDE solution
construction described above.
\end{abstract}
\begin{keywords}
\kwd{Stable L\'evy Diffusion}
\kwd{fractional diffusion}
\kwd{fractional advection-dispersion}
\kwd{heavy-tailed particle tracking}
\end{keywords}
%
%\begin{keywords}[MSC2010]%
%\kwd{}
%\kwd{}
%\kwd{}
%\end{keywords}

\received{\sday{15} \smonth{3} \syear{2018}}% Updated by
%VTEXPTS2LaTeX.exe, 14.06.2018 08:32
\revised{\sday{25} \smonth{5} \syear{2018}}% Updated by
%VTEXPTS2LaTeX.exe, 14.06.2018 08:32
\accepted{\sday{4} \smonth{6} \syear{2018}}% Updated by
%VTEXPTS2LaTeX.exe, 14.06.2018 08:32
\publishedonline{\sday{9} \smonth{7} \syear{2018}}
\end{frontmatter}\setcounter{footnote}{0}

%s1 #&#
\section{Introduction}
\label{s_intro}
\begingroup\abovedisplayskip=5pt\belowdisplayskip=5pt
The usual hydrological model for contamination/tracer transport through
a porous media is given by a second order advection-dispersion equation
(ADE) of the form $\frac{\partial c}{\partial t}=-\frac{\partial
}{\partial x}[v(x)c]+\frac{\partial^2}{\partial x^2}
[D(x)c ]$, where $c(x,t)$ is the tracer concentration at time
$t$, location $x$, $v$ is the drift velocity and $D$ is related to the
diffusivity of the media \cite{Bear,Fetter}. The probabilistic
approach to describe this flow from a mesoscopic view is given by the
hypothesis that the path of a randomly chosen tracer particle is a
Markov process that solves a stochastic differential equation (SDE)
driven by Brownian motion. The basis of this hypothesis is that the
conditional probability density function of this Markov process solves
a forward equation of the same form as the ADE \cite{Bhattacharya76,Bhattacharya83}.

However, for some heavy-tailed flows, the second order diffusion model
can be inadequate. For such cases a model called the fractional
advection-dispersion equation or fADE of the following form has been
proposed \cite{Zhang3}:
%
%e1 #&#
\begin{equation}
\label{eq:GeneralfADE} \dfrac{\partial c}{\partial t}=-\dfrac{\partial}{\partial
x}\bigl[v(x)c\bigr]+
\dfrac{\partial^{\alpha-1}}{\partial x^{\alpha-1}} \biggl[D(x)\dfrac{\partial c}{\partial x} \biggr].
\end{equation}
For the current discussion, we will consider a one-dimensional
concentration $c(x,t)$ where $x$ denotes the distance from the origin
of the plume, $v$ is the drift velocity and $D$ is a function that
changes with the diffusion of the tracer. The fractional
differentiation order $\alpha\in(1, 2)$ controls the tail of the
flow. Here the fractional derivative of order $\alpha$ for any
function $f$ is defined as in \cite{Baeumer} by:
\[
\frac{d^\alpha f(x)}{dx^\alpha}= \frac{-1}{\varGamma(1-\alpha)}\int \limits
^\infty_0
\bigl[f(x-u)-f(x)+uf'(x) \bigr]\alpha y^{-1-\alpha}du.
\]
The negative fractional derivative is given by $\frac{d^\alpha
f(x)}{d(-x)^\alpha}=\frac{d^\alpha g(x)}{dx^\alpha}$, where
$g(x)=f(-x)$.

If we assume that the diffusion coefficient $D$ is location invariant,
then the fADE can be associated with an SDE driven by an $\alpha
$-stable L\'evy process $X_t$ \footnote{Following the parametrization
in \cite{Taqqu} $X_t\sim \mathit{Stable}(\alpha,\beta,\mu,\sigma)$ where
$\alpha$ is the index of stability, $\beta$ is the skewness
parameter, $\mu=0$ is the location parameter and $\sigma=1$ is the
scale parameter.} of the form:
%
%e2 #&#
\begin{equation}
\label{eq:mainSD} dY_t=a(Y_t)dt+b(Y_t)dX_t.
\end{equation}
Following the Brownian diffusion argument, in the heavy-tailed
diffusion model, we assume a random particle's position at time $t$ is
given by the process $Y_t$ that solves SDE \eqref{eq:mainSD}. It can
be shown that $Y_t$ is a Markov process \cite{Applebaum}. Let us
denote the transition probability density function of $Y_t$ by
$p_{y_0}(y,t)$, i.e. $P(Y_t\in A|Y_0=y_0)=\int_Ap_{y_0}(y,t)dy$. For
the initial distribution $\mu(u)$ the pdf of $Y_t$ is given by
%
%e3 #&#
\begin{equation}
\label{eq:transitionpdf} p(y,t)=\int p_{u}(y,t)\mu(u)du,
\end{equation}
e.g., for a ground water tracer concentration modeling it is reasonable
to assume that all particles start at location $y_0$ at time $t=0$
hence $\mu(u)=I(u=y_0)$. If $1<\alpha<2$, it can be shown \cite
{Chakraborty} that $p(y,t)$ solves the forward equation:
%
%e4 #&#
\begin{align}
\label{forward3}
\frac{\partial p(y,t)}{\partial t}=-\frac{\partial}{\partial y} \bigl[a(y) p(y,t) \bigr]&+
\frac{(1+\beta)}{2} \biggl[\!-\cos \biggl(\frac
{\pi\alpha}{2} \biggr)\!
\biggr]^{-1}\!\!\cdot\frac{\partial^\alpha}{
\partial y^\alpha} \bigl[b^\alpha(y) p(y,t)
\bigr]
\nonumber
\\
&+\frac{(1-\beta)}{2} \biggl[\!-\cos \biggl(\frac{\pi\alpha}{2} \biggr)\!
\biggr]^{-1}\!\!\cdot\frac{\partial^\alpha} {\partial(- y)^\alpha
} \bigl[b^\alpha(y) p(y,t)
\bigr].
\end{align}
Typically the coefficient functions $a$ and $b$ depend on $t$ via $Y_t$
as in \eqref{eq:mainSD} hence the terms in the forward equation are
expressed as $a(y)$ and $b(y)$, where $y$ is supposed to be in the
range space of $Y_t$. However, without loss of generality, the
functions can also be considered as a mapping from $T\times\mathbb
{R}\rightarrow\mathbb{R}$ as $a(y,t)$ and $b(y,t)$ and the associated
forward equation will remain the same.

Thus for a location invariant $D$, the fADE in \eqref{eq:GeneralfADE}
can be written as:
%
%e5 #&#
\begin{equation}
\label{eq:fADE} \dfrac{\partial C}{\partial t}=-\dfrac{\partial}{\partial
x}\bigl[v(x,t)C\bigr]+D(t)
\dfrac{\partial^{\alpha}}{\partial x^{\alpha}} \biggl[\dfrac{\partial C}{\partial x} \biggr],
\end{equation}
which essentially has the same form as the forward equation given in
\eqref{forward3} with $\beta=1$, $a(x,t)=v(x,t)$, $ [-\cos
(\frac{\pi\alpha}{2} ) ]^{-1}b^\alpha(x,t)\equiv D(t)$.
This shows that in a heavy-tailed plume that follows fADE \eqref
{eq:fADE}, the position $Y_t$ of a randomly chosen tracer particle at
time $t$ solves SDE \eqref{eq:mainSD}. Choosing $\beta\neq1$ will
allow the particle transition in forward or backward direction.

In most practical tracking scenarios, the plume data consists of
observed tracer concentration values $c(x,t)$ over a range of locations
at certain time points. We write $p(y,t)=K_tc(y,t)$, where $K_t$ is a
suitable scale parameter that adjusts the total mass for $c(x,t)$ so
that $\int_{\mathbb{R}}p(y,t)dy=1$ for each $t$.

In case $v(\cdot)$ and $D(\cdot)$ are constants, \eqref{eq:mainSD}
reduces to $Y_t=at+bX_t$ with $a=v$, $b= [-\cos (\frac{\pi
\alpha}{2} )D ]^{1/\alpha}$ and consequently, $Y_t$
itself is a stable process. The parameters of this stable process then
can be estimated using the observed data (see \cite
{Chakraborty_Meerschaert} for details). However, when the coefficient
functions are not constants, the solution to the SDE in \eqref
{eq:mainSD} may not have a closed form.

In this paper, we present a numerical method to solve $Y_t$ through the
fitted fADE equation using the observed tracer concentration data.
Also, often for a groundwater contamination modeling problem the
functional form of the advection or dispersion coefficients are not
known to start with. In that case, the proper advection and dispersion
coefficient functions are needed to be estimated. This presents an
inverse problem. We formulate this inverse problem as an optimization
problem and develop a Levenberg--Marquardt (L-M) regularization method
to obtain the proper advection and dispersion coefficient function from
the observed data (the case for a linear drift function is presented
here). Then we proceed to generate $Y_t$ using these fitted coefficient
functions. The detailed methodology is described in Section~\ref{s_method}. Section~\ref{s:illustration}
provides an illustrative example. Section~\ref{sec4} includes concluding remarks.
\endgroup

%s2 #&#
\section{Methodology}
\label{s_method}
We propose a method for generating a weak solution for the SDE given in
\eqref{eq:mainSD}. The core technique is the standard probability
integral transformation. Suppose the probability density function (pdf)
of the solution process $Y_t$ at a given time point $t$ is given by
$p(y,t)$ and the cumulative distribution function (CDF) is defined by
$F_{Y_t}(x)=\int^x_{-\infty}p(y,t)dy$. Since $F(Y_t)\sim
\mathit{Uniform}(0,1)$ \cite{resnick2003probability}, a random observation of
$Y_t$, say $y$, can be generated by $F_{Y_t}^{-1}(u)=y$ where $u$ is a
random observation from $\mathit{Uniform}(0,1)$. For the present problem, closed
form of $p(t,y)$ is not available, and we use point-wise numerical
approximations to estimate $p(y,t)$ for any fixed $t$ and fixed $y$.
Therefore we can only use a numerical scheme to approximate $F_{Y_t}$
and $F_{Y_t}^{-1}$.

%s2.1 #&#
\subsection{Simulating a solution process for an SDE driven by a
heavy-tailed stable process}
\label{subs:simulation}
For a heavy-tailed plume where the particle path can be modeled by the
SDE in \eqref{eq:mainSD}, the forward equation is of the form \eqref
{eq:fADE}. To generate the solution to this SDE we can consider two
scenarios: (i) we may want to solve %for
an SDE with a given form of coefficient functions $v(\cdot)$ and
$D(\cdot)$ and driven by a given $\alpha$-stable process; (ii) a more
practical problem when only the tracer concentrations at different time
points over a range of locations are observed from a plume, while the
coefficient functions need to be estimated along with
the parameters of the driving process $X_t$ in \eqref{eq:mainSD} from
the observed concentration values.

In either case, first, we find $p(t,y)$ that solves the forward
equation \eqref{forward3} with given or estimated coefficient
functions and then generate $Y_t$ using these $p(t,y)$'s.

\smallskip

\textbf{\underline{$Y_t$ Simulation Steps}:}
\begin{description}
\item[Step 1: Numerical Solution for $p(t,y)$:]
Assuming the probability density $p(y,t)=K_tc(y,t)$, we estimate $K_t$,
$\alpha$ and $\beta$ from the observed data that assumed to satisfy
fADE \eqref{eq:fADE}. $p(y,t)$ can be solved numerically using a
discretization scheme from the forward equation \eqref{forward3}. In
case the coefficient functions in \eqref{forward3} are not known, an
inverse problem optimization technique can be used for estimating the
parameters. See Section~\ref{ssNumForwd} for details of this step.

In either case the discretization method provides only point-wise
approximation of $p(y,t)$ values for given $t$ and $y$ and not a closed
functional form.

%\item[]\textbf{Step 2: Generating $Y_t$:}
\item[Step 2: Generating $Y_t$:]
To generate $Y_t$ for a fixed $t$ using an inverse CDF transformation,
we present a straightforward scheme with a trapezoidal rule to
approximate $F_{Y_t}$ and $F_{Y_t}^{-1}$.
Approximate $Y_t$'s can be simulated using $p(y,t)$'s from step 1 as follows:
\begin{enumerate}
\item[(a)] Choose suitable grid $y_0< y_1<\cdots< y_n$ and solve for
$p(y_i,t)$ as mentioned in step 1. The endpoints $y_0$ and $y_n$ can be
chosen so that $p(y,t)$ is small or negligible for $y<y_0$ or $y>y_n$.
The fitted coefficient and parameters from the observed data in step 1
can be used for any $y_i$'s.
\item[(b)] Approximate $F_{Y_t}(y_i)$ by $\hat{F}_{Y_t}(y_i)$ for
$i=1,2,\ldots, n$, using the trapezoidal rule:
\begin{align*}
\hat{F}_{Y_t}(y_1)&=(y_1-y_0)
\dfrac{p(y_1,t)+p(y_0,t)}{2},
\\
\hat{F}_{Y_t}(y_2)&= \hat{F}_{Y_t}(y_1)+(y_2-y_1)
\dfrac
{p(y_2,t)+p(y_1,t)}{2},
\\
\vdots&
\\
\hat{F}_{Y_t}(y_r)&= \hat{F}_{Y_t}(y_{r-1})+(y_r-y_{r-1})
\dfrac
{p(y_r,t)+p(y_{r-1},t)}{2},
\\
\vdots&
\\
\hat{F}_{Y_t}(y_n)&= \hat{F}_{Y_t}(y_{n-1})+(y_n-y_{n-1})
\dfrac
{p(y_n,t)+p(y_{n-1},t)}{2}.
\end{align*}
\item[(c)] For a randomly generated $\mathit{Uniform}(0,1)$ observation $u$,
find its placement in the grid, i.e. $y_r$, $y_{r+1}$ such that $ \hat
{F}_{Y_t}(y_r)\leq u< \hat{F}_{Y_t}(y_{r+1})$, for some $r\in\{
0,1,\ldots,n\}$. A way of generating associated random observation of
$Y_t$ or $F^{-1}_{Y_t}(u)$ (approximated) as $\hat{y}$ is described below:
\[
\hat{y}=\lleft\{ %
\begin{array}{ll}
y_r, & \hbox{if }\ u=y_r; \\
y_r+\dfrac{u-\hat{F}_{Y_t}(y_r)}{p(t,y_r)}, & \hbox{if }\ y_r< u<y_{r+1}.
\end{array} %
\rright.
\]
In case $u<\hat{F}_{Y_t}(y_0)$ set $\hat{y}=y_0$ and if $u\geq\hat
{F}_{Y_t}(y_{n})$ set $y=y_{n}$.
\item[(d)] To generate $N$ i.i.d approximated random observations
$\hat{y}_1,\hat{y}_2,\ldots, \hat{y}_N$ of $Y_t$, start with
randomly generated independent $\mathit{Uniform}(0,1)$ observations $u_1,
u_2,\ldots u_N$ and then repeat step (c). Note that the efficiency of
this approximation will depend on the spacing of the chosen grid since
we are actually generating sample from distribution $\hat{F}_{Y_t}$
here and $\hat{F}_{Y_t}\rightarrow F_{Y_t}$ as
$\max\limits_{0\leq i\leq n}|y_{i+1}-y_i|\rightarrow0$ following the
convergence property of the standard trapezoidal rule.
\end{enumerate}
\end{description}

%s2.2 #&#
\subsection{A formulation for a numerical solution of the forward equation}
\label{ssNumForwd}
In this subsection, we will present a finite volume scheme with the
characteristic line for solving a forward equation of the form \eqref
{forward3}. To incorporate a discretized finite volume scheme with
proper boundary condition, we re-parameterize the forward equation with
fractional order $2-\lambda$ with $0<\lambda<1$ (consistent with the
construction described in \cite{Ref:J:WangH+AlLawatiaM:nmpde:y:2006})
as follows:
%
%e6 #&#
\begin{equation}
\label{eq:Xu1} %
\begin{aligned}
&\f{\p p} {\p t} + \f{\p} {\p x} a(x)p - \frac{\p}{\p x} \bigl[ b(t) \bigl(\gamma{}_{x_l}D_{x}^{-\lambda}+(1-\gamma ){}_{x}D_{x_r}^{-\lambda}\bigr)  Dp \bigr] = 0, \quad x_l < x < x_r,\\
&p(x_l,t) =0, \qquad p(x_r,t) = 0, \quad 0 < t < T.
\end{aligned} %
\end{equation}
Essentially, compared to \eqref{forward3} the fractional derivative
order $\alpha=2-\lambda$, while $\gamma=\frac{1+\beta}{2}$
indicates the relative weight of forward versus backward transition.
${}_{x_l}D_{x}^{-\lambda}$ and ${}_{x}D_{x_r}^{-\lambda}$ represent
the left and right fractional integral operators defined as
%
%e7 #&#
\begin{equation}
\label{eq:Xu2} %
\begin{aligned} {}_{x_l}D_{x}^{-\lambda}
f(x):&=\frac{1}{\varGamma(\beta)}\int_{x_l}^{x}
(x-s)^{\lambda-1} f(s){\di}s,\\%[0.15in]
{}_{x}D_{x_r}^{-\lambda} f(x):&=
\frac{1}{\varGamma(\beta)}\int_{x}^{x_r} (s-x)^{\lambda-1}
f(s){\di}s. \end{aligned} %
\end{equation}
To comply with the boundary conditions, $x=x_l$ and $x=x_r$ are set to
be the inflow and outflow boundaries respectively \cite
{Ref:J:WangH+AlLawatiaM:nmpde:y:2006}, while $a(x_l)$ and $a(x_r)$ are
assumed to be non-negative.

%s2.2.1 #&#
\subsubsection{The discretized finite volume scheme and the
accumulation terms}

Let us define a partition $\pt$ on the time interval $[0,T]$, and a
partition $\px$ on the space interval $[x_l,x_r]$
%
%e8 #&#
\begin{equation}
\label{eq:Xu3} %
\begin{aligned}
\pt: 0      &= t_0 < t_1 < \cdots< t_N = T,\\
\px: x_l    &= x_0 < x_1 < \cdots<x_I =x_r,
\end{aligned} %
\end{equation}
with $\Delta t_n = t_n - t_{n-1}$ for $n=1,\ldots,N$ and $\Delta x_i=
x_i - x_{i-1}$ for $i=1,\ldots,I$. Let $\bar\px$ be the dual
partition of $\px$ defined by $\xim= \f1{2}(x_{i-1}+x_i)$ with
$\Delta\xip= \xip- \xim$ for $i=1,\ldots,I-1$ as well as $\Delta
x_{1/2} =x_{1/2} - x_0$ and $\Delta x_{I+1/2}=x_I - x_{I-1/2}$.

Let $y=r(\theta;\bar x,\bar t)$ be a continuous and piecewise-smooth
curve that passes through the point $\bar x$ at time $\bar t$, such
that $r(\theta;\xim,t_n)$ and $r(\theta;\xip,t_n)$ do not meet each
other during the time period $[t_{n-1},t_n]$. We define a space-time
control volume $\omin$ by extending the cell $[\xim,\xip]$ along the
curve $r(t;x,t_n)$ from $t=t_n$ to $t=t_{n-1}$
%
%e9 #&#
\begin{equation}
\label{eq:Xu4} \omin=\bigl\{(x,t): r(t;\xim,t_n)< x< r(t;
\xip,t_n),\quad t_{n-1}< t < t_n\bigr\}.
\end{equation}
Assuming the prism $\omin$ does not intersect the boundary $x=x_l$ and
$x=x_r$ of the domain during the time period $[t_{n-1},t_n]$, let
$x^*=r(t_{n-1};x,t_n)$ be the foot at time $t_{n-1}$ of the curve with
head $x$ at time $t_n$.

Integrating equation \eqref{eq:Xu1} over the control volume $\omin$
we get
%
%e10 #&#
\begin{equation}
\label{eq:Xu5} \int_{\omin} \f{\p p} {\p t} dxdt +\int
_{t_{n-1}}^{t_n} \Bigl(a(x) p - b(t)\f{\p^{1-\lambda}} {
\p x^{1-\lambda}} p \Bigr) \Big|_{r(t;\xim,t_n)}^{r(t;\xip,t_n)} dt = 0.
\end{equation}
Without loss of generality, the accumulation term in \eqref{eq:Xu4}
can be evaluated by assuming
%
%e11 #&#
\begin{equation}
\label{eq:Xu6} \xim^* < \xim< \xip^* < \xip,
\end{equation}
and accumulation term can be re-written as follows:
%
%e12 #&#
\begin{align}
\label{eq:Xu7} \int_{\omin} \f{\p p} {\p t} dxdt & =\int
_{\xim^*}^{\xim} \Biggl[ \int_{t_{n-1}}^{t(x;\xim,t_n)}
\f{\p p} {\p t} \:dt \Biggr]\: dx\nonumber\\%[0.1in]
&  \quad +\int_{\xim}^{\xip^*} \Biggl[ \int
_{t_{n-1}}^{t_n} \f{\p p} {\p t} \: dt \Biggr] \:dx + \int
_{\xip^*}^{\xip} \Biggl[ \int_{t(x;\xip,t_n)}^{t_n}
\f{\p p} {\p t} \: dt \Biggr]\: dx.
\end{align}
Here the notation $t(x;\xim,t_n)$ represents the time instant that the
curve\hfill\break
$r(t;,\xim,t_n) = x$. The notation $t(x;\xip,t_n)$ is defined
similarly.

A simple calculation of the second term on the right-hand side of
\eqref{eq:Xu7} yields
%
%e13 #&#
\begin{equation}
\label{eq:Xu8} \int_{\xim}^{\xip^*} \Biggl[ \int
_{t_{n-1}}^{t_n} \f{\p p} {\p t} \: dt \Biggr]dx = \int
_{\xim}^{\xip^*}p(x,t_n)dx-\int
_{\xim}^{\xip^*} p(x,t_{n-1})dx.
\end{equation}
The first and third terms on the right-hand side of \eqref{eq:Xu7} are
integrated as follows
%
%e14 #&#
\begin{align}
\label{eq:Xu9} &\int_{\xim^*}^{\xim} \Biggl[ \int
_{t_{n-1}}^{t(x;\xim,t_n)} \f{\p p} {\p t} \: dt \Biggr]dx\nonumber\\%[0.1in]
&\quad =\int_{t_{n-1}}^{t_n} p \bigl(r(t;
\xim,t_n),t \bigr) \f{\p r(t;\xim,t_n)} {\p t}\: dt - \int
_{\xim^*}^{\xim} p(x,t_{n-1}) \: dx,
\end{align}
and
%
%e15 #&#
\begin{align}
\label{eq:Xu10} &\int_{\xip^*}^{\xip} \Biggl[ \int
_{t(x;\xip,t_n)}^{t_n} \f{\p p} {\p t} \: dx \Biggr]dt\nonumber\\%[0.2in]
&\quad = \int_{\xip^*}^{\xip} p(x,t_n) \:
dx- \int_{t_{n-1}}^{t_n} p \bigl(r(t;\xip,t_n),t
\bigr)\f{\p r(t;\xip,t_n)} {\p t}~dt.
\end{align}

Incorporating equations \eqref{eq:Xu8}--\eqref{eq:Xu18} into \eqref
{eq:Xu7} we get
%
%e16 #&#
\begin{equation}
\label{eq:Xu11} %
\begin{aligned} \int_{\omin} \f{\p
p} {\p t} dxdt & =\int_{\xim}^{\xip}p(x,t_n)dx-
\int_{\xim^*}^{\xip^*} p(x,t_{n-1})dx\\%[0.2in]
&\quad + \int_{t_{n-1}}^{t_n} p \bigl(r(t;
\xim,t_n),t \bigr) \f{\p r(t;\xim,t_n)} {\p t} \: dt\\%[0.15in]
& \quad - \int_{t_{n-1}}^{t_n} p \bigl(r(t;
\xip,t_n),t \bigr) \f{\p r(t;\xip,t_n)} {\p t} \: dt.
\end{aligned} %
\end{equation}
The derivation shows that \eqref{eq:Xu11} does not depend on the
assumption \eqref{eq:Xu4}.

Now we just set the space-time boundaries $r(t;\xpm,t_n)$ of the
control volume $\omin$ to be the characteristic curves, which are
defined by the initial value problem of the ordinary differential equation
%
%e17 #&#
\begin{equation}
\label{eq:Xu12} \f{d r} {dt} = a(r,t), \qquad r(t;x,t_n) ~
|_{t=t_n} = x.
\end{equation}
That is, because the characteristics $r(t;\xpm,t_n)$ are assumed to be
tracked exactly and hence the residual advection term vanishes
naturally.

Substituting \eqref{eq:Xu11} for the accumulation term in \eqref
{eq:Xu5}, we obtain a locally conservative reference equation on an
interior space-time control volume $\omin$ as follows:
%
%e18 #&#
\begin{equation}
\label{eq:Xu13} %
\begin{aligned}
&\int_{\xim}^{\xip}p(x,t_n)
\:dx + \int_{t_{n-1}}^{t_n} F(p,b) \bigl(r(t;
\xip,t_n),t \bigr) \: dt\\%[0.1in]
&\qquad - \int_{t_{n-1}}^{t_n} F(p,b) \bigl(r(t;
\xim,t_n),t \bigr) \: dt\\%[0.15in]
&\quad = \int_{\xim^*}^{\xip^*}p(x,t_{n-1})~dx,
\end{aligned} %
\end{equation}
which can be approximated as
%
%e19 #&#
\begin{equation}
\label{eq:Xu14} %
\begin{aligned}
&\int_{\xim}^{\xip}p(x,t_n)\:dx + \Delta t \bigl[F(p,b) (\xip,t_n )-F(p,b) (\xim,t_n ) \bigr]\\
&\quad = \int_{\xim^*}^{\xip^*}p(x,t_{n-1})~dx,
\end{aligned} %
\end{equation}
with the Lagrangian interface fluxes
%
%e20 #&#
\begin{equation}
\label{eq:Xu15}  F(p,b) \bigl(r(t;\xpm,t_n),t \bigr) := -b~\f{
\p^{1-\lambda} } {\p x^{1-\lambda}} p \bigl(r(t;\xpm ,t_n),t
\bigr).
\end{equation}

Note that the interface fluxes $F(p,b)(r(t;\xpm,t_n),t)$ are defined
across the space-time boundary $r(t;\xim,t_n)$ and $r(t;\xip,t_n)$ of
the space-time control volume $\omin$.

To discretize equation \eqref{eq:Xu1}, we further let $\{\phi_i\}
_{i=1}^{I}$ be a set of hat functions such that $\phi_i(x_i)=1$ and
$\phi_i(x_j)=0$ for $j\neq i$.
Hence the finite volume approximation $p_h$ to the true solution $p$
can be expressed as
\[
p_h(x) = \sum_{j=1}^{I-1}
p_j\phi_j(x). %+ p_l\phi_0(x) + p_r\phi_I(x),
\]
and the finite volume scheme on the interval $[x_{i-1/2}, x_{i+1/2}]$
for $i=1,\ldots,I-1$ has the form
%
%e21 #&#
\begin{equation}
\label{eq:Xu16}  \int_{\xim}^{\xip}p(x,t_n)
\:dx + \Delta t \sum_{j=1}^{I-1}
p_{j}^{n} z_{i,j} = \int_{\xim^*}^{\xip^*}
p(x,t_{n-1})~dx,
\end{equation}
where, $[z_{i,j}]_{i,j=1}^{I-1}$ is the coefficient matrix of the
fractional term, and is given by
%
%e22 #&#
\begin{equation}
\label{eq:Xu17} %
\begin{aligned} z_{i,j} = b \bigl[ \bigl(
\gamma{}_{0}D_{x_{i-1/2}}^{-\lambda
}+(1-\gamma){}_{x_{1-1/2}}D_{1}^{-\lambda}
\bigr)  D\phi_j(x_{i-1/2,
t_n})\\%[0.15in]
\qquad\quad - \bigl(\gamma{}_{0}D_{x_{+1/2}}^{-\lambda}+(1-
\gamma ){}_{x_{1+1/2}}D_{1}^{-\lambda}\bigr)  D
\phi_j(x_{i+1/2, t_n}) \bigr]. \end{aligned} %
\end{equation}

Assuming the trial functions $p(x,t_n)$ are chosen to be piecewise
linear functions on $[x_l,x_r]$ with respect to the fixed spatial
partition $\px$, we evaluate the accumulation term at time step $t_n$
in the finite volume scheme \eqref{eq:Xu16} analytically as
%
%e23 #&#
\begin{equation}
\label{eq:Xu18} \int_{\xim}^{\xip}p(x,t_n)
\:dx = \f1{8} \bigl[ \Delta x_i (\pim+3\ppi)+ \Delta
x_{i+1}(3\ppi+ \pip) \bigr].
\end{equation}

In addition, we can compute $\xim^*$ and $\xip^*$ at time step
$t_{n-1}$ via a backward approximate characteristic tracking. Since the
trial function $p(x,t_{n-1})$ is also piecewise linear with respect to
the fixed spatial partition $\px$ at time $t_{n-1}$, we can evaluate
the accumulation term at time step $t_{n-1}$ analytically. Because the
accumulation term at time step $t_{n-1}$ affects only the right-hand
side of the finite volume scheme \eqref{eq:Xu16}, the scheme retains a
symmetric and positive-definite coefficient matrix. Furthermore, the
finite volume scheme \eqref{eq:Xu16} is locally conservative, even if
the characteristics are computed approximately.

%s2.2.2 #&#
\subsubsection{The fractional diffusion term and the stiffness matrix}

The finite volume scheme~(\ref{eq:Xu16}) appears similar to the one for the
canonical second-order diffusion equation, but has a fundamental
difference. Although the hat functions $\phi_j$ have local support,
\eqref{eq:Xu2} reveals that ${}_{x_l}D_{x}^{-\lambda} \phi_j$ and
${}_{x}D_{x_r}^{-\lambda} \phi_j$ have global support. Therefore, the
stiffness matrix $Z:=[z_{i,j}]_{i,j=1}^{I-1}$ is a full matrix, which
requires $O(N^2)$ of memory to store. Numerical schemes for
space-fractional differential equations were traditionally solved by
the Gaussian type direct solvers that require $O(N^3)$ of computations
\cite{Ref:J:NegreteD+CarrerasBA+LynchVE:pp:y:2004,Ref:J:ErvinVJ+RoopJP:nmpde:y:2005,Ref:J:LiuF+AnhV+TurnerI:jcam:y:2004,Ref:J:MeerschaertM+TadjeranC:jcam:y:2004}. In recent years, there have
been some other notable developments in methods for solving the
algebraic linear systems arising from discretization of
fractional-order problems, especially for space dimension higher than
one (see \cite{Ref:Breiten}); we plan to explore them in our future
work.

To simplify the computation of the diffusion term, we have to explore
the structure of the stiffness matrix $Z$ \cite
{Ref:J:JiaJ+WangH:jcp:y:2015,Ref:J:WangH+DuN:jcp:y:2013}. The entries
of the stiffness matrix $Z$ are presented below.

Its diagonal entries are given by:
%
%e24 #&#
\begin{equation}
\label{eq:Xu19}
\begin{aligned}
z_{i,i}& =\f{1}{\varGamma(\lambda+1)2^\lambda} b\gamma h^{\lambda-1}+\f{1}{\varGamma(\lambda+1)} b(1-\gamma) h^{\lambda-1}\bigg(2\left(\f{1}{2}\right)^\lambda-\left(\f{3}{2}\right)^\lambda\bigg)\\%[0.05in]
& \quad +\f{1}{\varGamma(\lambda+1)} b\gamma h^{\lambda-1}\bigg(2\left(\f{1}{2}\right)^\lambda-\left(\f{3}{2}\right)^\lambda\bigg)\\%[0.05in]
& \quad +\f{1}{\varGamma(\lambda+1)2^\lambda} b(1-\gamma) h^{\lambda-1}\\%[0.15in]
&= \f{1}{\varGamma(\lambda+1)} b h^{\lambda-1} \bigg(2^{-\lambda}+2\left(\f{1}{2}\right)^\lambda-\left(\f{3}{2}\right)^\lambda\bigg),
\end{aligned}
\end{equation}
for $1\leq i \leq I-1$.

Whereas, the sub-triangular entries of the matrix $Z$ are given by:
%
%e25 #&#
\begin{equation}
\label{eq:Xu20}
\begin{aligned}
z_{i,i-1}& = -\f{1}{\varGamma(\lambda+1)} b\gamma h^{\lambda-1}\bigg(2\left(\f{1}{2}\right)^\lambda-\left(\f{3}{2}\right)^\lambda\bigg)-\f{1}{\varGamma(\lambda+1)2^\lambda} b(1-\gamma) h^{\lambda-1}\\%[0.15in]
&\quad -\f{1}{\varGamma(\lambda+1)} b\gamma h^{\lambda-1}\bigg(\left(\f{5}{2}\right)^\lambda-2\left(\f{3}{2}\right)^\lambda+\left(\f{1}{2}\right)^\lambda\bigg)\\%[0.15in]
&= \f{1}{\varGamma(\lambda+1)} b h^{\lambda-1}\bigg(3\left(\f{3}{2}\right)^\lambda\gamma-3\left(\f{1}{2}\right)^\lambda\gamma-\left(\f{5}{2}\right)^\lambda\gamma-2^{-\lambda}(1-\gamma)\bigg),
\end{aligned}
\end{equation}
for $2\leq i \leq I-1$, and
%
%e26 #&#
\begin{equation}
\label{eq:Xu21}
\begin{aligned}
z_{i,j} &= \f{1}{\varGamma(\lambda+1)} b\gamma h^{\lambda-1} \bigg(\left(i-j+\f{1}{2}\right)^\lambda-2\left(i-j-\f{1}{2}\right)^\lambda+\left(i-j-\f{3}{2}\right)^\lambda\bigg)\\%[0.15in]
& \quad -\f{1}{\varGamma(\lambda+1)} b\gamma h^{\lambda-1} \bigg(\left(i-j+\f{3}{2}\right)^\lambda-2\left(i-j+\f{1}{2}\right)^\lambda+\left(i-j-\f{1}{2}\right)^\lambda\bigg)
\end{aligned}
\end{equation}
for $3\leq i \leq I-1$ and $1\leq j \leq i-2$.

The super-triangular entries of matrix $Z$ can be also derived as\vadjust{\vfill{\eject}}
%
%e27 #&#
\begin{align}
\label{eq:Xu22}
z_{i,i+1} &= -\f{1}{\varGamma(\lambda+1)} b(1-\gamma) h^{\lambda-1}\bigg(\left(\f{1}{2}\right)^\lambda-2\left(\f{3}{2}\right)^\lambda+\left(\f{5}{2}\right)^\lambda\bigg)\nonumber\\%[0.15in]
&\quad -\f{1}{\varGamma(\lambda+1)2^\lambda} b\gamma h^{\lambda-1} -\f{1}{\varGamma(\lambda+1)} b(1-\gamma) h^{\lambda-1} \bigg(2\left(\f{1}{2}\right)^\lambda-\left(\f{3}{2}\right)^\lambda\bigg)\nonumber\\%[0.15in]
&= \f{1}{\varGamma(\lambda+1)} b h^{\lambda-1}\bigg(3\left(\f{3}{2}\right)^\lambda(1-\gamma)-3\left(\f{1}{2}\right)^\lambda(1-\gamma)\nonumber\\
&\quad -\left(\f{5}{2}\right)^\lambda(1-\gamma)-2^{-\lambda}\gamma\bigg),
\end{align}
for $1\leq i \leq I-2$, and
%
%e28 #&#
\begin{align}
\label{eq:Xu23}
z_{i,j} &= \f{1}{\varGamma(\lambda+1)} b(1-\gamma) h^{\lambda-1}\bigg(2\left(j-i+\f{1}{2}\right)^\lambda\nonumber\\
&\quad -\left(j-i-\f{1}{2}\right)^\lambda-\left(j-i+\f{3}{2}\right)^\lambda\bigg)\nonumber\\%[0.15in]
&\quad -\f{1}{\varGamma(\lambda+1)} b(1-\gamma) h^{\lambda-1}\bigg(2\left(j-i-\f{1}{2}\right)^\lambda\nonumber\\
&\quad -\left(j-i-\f{3}{2}\right)^\lambda-\left(j-i+\f{1}{2}\right)^\lambda\bigg)
\end{align}
for $1\leq j \leq I-3$ and $j+2\leq i \leq I-1$.

%s2.2.3 #&#
\subsubsection{Estimation of coefficient function}

We present the case where the drift coefficient function in fADE \eqref
{eq:fADE} (and therefore in equation \eqref{eq:Xu1}) is assumed to be
a piecewise linear function of $x$,
%
%e29 #&#
\begin{equation}
\label{eq:Xu24} a(x) = \lleft\{ %
\begin{array}{ll}
a_0 - a_1 x, & \textrm{if $x_l\leq x\leq x_m$},\\
a_2 - a_3 x, & \textrm{if $x_m< x\leq x_r$},
\end{array} %
\rright.
\end{equation}
where parameters $a_0$, $a_1$, $a_2$ and $a_3$ are to be estimated from
the observed concentration data. The main idea on this part is to
obtain certain measurements through physical or mechanical experiments,
and then use the data to calibrate these parameters in the fADE (see
\cite{Ref:J:FuH+WangH+WangZ:jsc:y:2018}). This is an inverse problem:
based on the initial guess $p_0$ of the equation \eqref{eq:Xu1}, and
certain observation (or desired) data such as values of the state
variable $\mathbf{g}$ at the final time, we attempt to seek for the
constant parameters $a_0$, $a_1$, $a_2$ and $a_3$ from the governing
differential equation \eqref{eq:Xu1}.

We formulate the inverse problem as an optimization and develop a
Levenberg--Marquardt (L-M) regularization method (see, \cite
{Ref:B:ChaventG:sN:y:2009,Ref:B:NocedalJ+WrightSJ:sNY:y:2006,Ref:B:SunW+YuanY:sNY:y:2006}) to
iteratively identify the parameter. It is known that the inverse
problem usually requires multiple runs of the forward problem.
Considering the computational cost of the forward problem is already
high, the inverse problem could become infeasible. Hence we propose an
optimization algorithm for the parameter estimation.

Here we only present the details of the fitting of $a_0$ and $a_1$. The
other two parameters $a_2$ and $a_3$ can be estimated similarly. The
parameter identification of $a_0$ and $a_1$ can be formulated as follows:
let $\balpha:= \{a_0,a_1\}$, then to find $\alpha_{inv}$ that satisfies
%
%e30 #&#
\begin{equation}
\label{eq:Xu25}  \alpha_{inv} = \text{arg}~~\min_{\balpha}~~
\gm(\balpha) := \f{1} {2}\sum_{i=1}^{M}
\bigl[p(x_i,t;\balpha)-\hat{p}(x_i,t)
\bigr]^2,
\end{equation}
where, $\hat{p}(x_i,t)=\hat{K}_tc(x_i,t)$ for observed concentration
$c$ at location $x_i$, time $t$.

In case the data is available for time points $t_1, t_2,\ldots, t_R$
we rewrite equation \eqref{eq:Xu25} as
%
%e31 #&#
\begin{equation}
\label{eq:Xu26} %
\begin{aligned} \alpha_{inv} &=
\text{arg}~~\min_{\balpha}~~\gm(\balpha)
\\
&= \f{1} {2}\sum_{k=1}^{R}\sum
_{i=1}^{M_k} w_k \bigl[p(x_{t_k,i},t_k;
\balpha)-\hat{p}(x_{t_ki},t_k) \bigr]^2,
% \\
%&+ \f{1}{2}\sum_{i=1}^{M_2} w_2 \Big[p(x_{2,i},t_2;\balpha)-
%\hat{p}(x_{2,i},t_2)\Big]^2 \\
% &+ \f{1}{2}\sum_{i=1}^{M_3} w_3 \Big[p(x_{3,i},t_3;\balpha)-
%\hat{p}(x_{3,i},t_3)\Big]^2,
\end{aligned} %
\end{equation}
here $x_{t_k,i}$ is the $i${th} observed location at time $t_k$ and
$w_k$ is the weight assigned for the sample set available for the same
time point. (For example the data set used in Section~\ref
{s:illustration} has $R=3$).

An iteration algorithm such as the Newton method with line searching
could be employed to find the solution to \eqref{eq:Xu26}. Basically,
the Newton algorithm for minimizing \eqref{eq:Xu26} uses the first and
second derivatives of the objective function $\gm(\balpha)$:
%
%e32 #&#
\begin{equation}
\label{eq:Xu27}  \balpha_{k+1} = \balpha_{k} - \f{
\gm'(\balpha_{k})} {\gm''(
\balpha_{k})},
\end{equation}
where $k$ represents the $k$th iteration.
It is easy to check that \eqref{eq:Xu27} is equivalent to solve
%
%e33 #&#
\begin{equation}
\label{eq:Xu28}  \balpha_{k+1} = \balpha_{k} - \bigl(
\mathbf{J}_{k}^{T} \mathbf{J}_{k}
\bigr)^{-1} \mathbf{J}_{k}^{T}
\mathbf{r}_{k},
\end{equation}
where
%
%e34 #&#
\begin{equation}
\label{eq:Xu29}
\begin{aligned}
\mathbf{J}_{k}
&= \big(\mathbf{J}^{1};\mathbf{J}^{2}; \ldots;\mathbf{J}^{R}\big
)^T,\quad\text{with}\\%[0.15in]
\mathbf{J}^{i}
&= \left( \f{\p p(x_{t_i,1},t_i;\balpha)}{\p\balpha},\ldots,\f
{\p p(x_{t_i,M_i},t_i;\balpha)}{\p\balpha} \right),
\end{aligned}
\end{equation}
and
%
%e35 #&#
\begin{equation}
\label{eq:Xu30} %
\begin{aligned}
\mathbf{r}_{k}&=\bigl(\mathbf{r}^{1};\mathbf{r}^{2};\ldots;\mathbf{r}^{R}\bigr)^T, \quad \text{with}\\%[0.15in]
\mathbf{r}^{i}&= (r_{i,1},\ldots,r_{i,M_i}),\quad  r_{i,j}=p(x_{t_i,j},t_i;\balpha)-
\hat{p}(x_{t_i,j},t_i), \end{aligned} %
\end{equation}
for $j=1,\ldots,M_i$, and $i=1,2,\ldots, R$. Note that in practice,
we always use the finite difference
\[
 \f{p(x_{t_i,j},t_i;\balpha+\delta)-p(x_{t_i,j},t_i;
\balpha)} {\delta}
\]
with a small enough $\delta$ to approximate the derivatives in \eqref
{eq:Xu29}.

However, the Newton method may fail to work because of $\mathbf
{J}_{k}^{T} \mathbf{J}_{k}$ may be nearly zero. Therefore, the search direction
$d_k:=-\mathbf{J}_{k}^{T} \mathbf{r}_{k}/\mathbf{J}_{k}^{T} \mathbf
{J}_{k}$ may not point in a decent direction.

A common technique to overcome this kind of problem is the L-M
algorithm (or Levenberg algorithm since a single parameter case is
considered in this paper), which modifies \eqref{eq:Xu27} by the
following formulation
%
%e36 #&#
\begin{equation}
\label{eq:Xu31}  \balpha_{k+1} = \balpha_{k} - \bigl(
\mathbf{J}_{k}^{T} \mathbf{J}_{k} + \varrho
_k I_2 \bigr)^{-1} \mathbf{J}_{k}^{T}
\mathbf{r}_{k},
\end{equation}
where $\varrho_k$ is a positive penalty parameter, and $I_2$ is a
$2\times2$ identity matrix. The method coincides with the Newton
algorithm when $\varrho_k=0$; and it gives a step closer to the
gradient descent direction when $\varrho_k$ is large.

\medskip

\textbf{\underline{Algorithm A: Parameter Identification Algorithm}:}

\medskip

Given the observation data $\mathbf{g}$ and an initial guess $\balpha
_0=\{a_0,a_1\}_0$, choose $\rho\in(0,1)$, $\sigma\in(0,1/2)$,
$\varrho_0>0$, and $\delta$ small enough.
For $k=0,1,\ldots,M$, (let $M=M_1+M_2+M_3$),
\begin{enumerate}
\item[(1)] Solve the equation \eqref{eq:Xu1} corresponding to
$\balpha_{k}$ and $\balpha_{k}+\delta$ respectively to obtain
$p(\cdot,\cdot,\balpha_{k})$ and $p(\cdot,\cdot,\balpha
_{k}+\delta)$.
\vspace{1mm}
\item[(2)] Compute $\mathbf{J}_{k}$ and $\mathbf{r}_{k}$, and update
the search direction
\[
 \mathbf{d}_k:=- \bigl(\mathbf{J}_{k}^{T}
\mathbf{J}_{k} + \varrho_k I_2
\bigr)^{-1} \mathbf{J}_{k}^{T}
\mathbf{r}_{k}.
\]
\item[(3)] Determine the search step $\rho^m$ by Armijo rule:
\[
 \gm \bigl(\balpha_{k}+\rho^m \mathbf{d}_k
\bigr) \leq\gm(\balpha_{k}) + \sigma\rho^m
\mathbf{d}_k \mathbf {J}_{k}^{T}
\mathbf{r}_{k}
\]
where $m$ is the smallest non-negative integer.
\vspace{1mm}
\item[(4)] If $|\rho^m \mathbf{d}_k|\leq\text{Tolerance level}$,
then stop and let $\alpha_{inv}:=\balpha_k$.
Otherwise, update
\[
 \balpha_{k+1} := \balpha_{k} +\rho^m
\mathbf{d}_k, \qquad \varrho_{k+1}:=\varrho_k
/2,
\]
and go to the first step for the iteration again.
\end{enumerate}

Algorithm A summarizes the proposed parameter estimation steps through
the inverse problem approach which includes the details of the L-M
method. In particular, the Armijo rule \cite{Ref:J:ArmijoL:pjm:y:1966}
in the third step, known as one of the inexact line search techniques,
is imposed to ensure the objective function $\gm$ has sufficient
descent. Other rules and the related convergence theory can be found in
\cite{Ref:B:SunW+YuanY:sNY:y:2006}.

%s3 #&#
\section{An illustrative example}
\label{s:illustration}
In this section, we use a heavy-tailed groundwater tracer concentration
data to illustrate the methodology described in Section~\ref
{s_method}. The data comes from natural-gradient tracer tests conducted
at the MacroDispersion Experimental (MADE) site at Columbus Air Force
Base in northeastern Mississippi, precisely the MADE-2 tritium plume
data \cite{Boggs93}. %The data has four spatial snapshots at day 27,
%day 132, day 224, and day 328 days after injection.
The data consists of the maximum concentration measured in vertical
slices perpendicular to the direction of plume travel, at day 27, day
132, day 224, and day 328 days after injection (see Figure 4 in \cite
{Benson2001}).

In \cite{Chakraborty_Meerschaert}, an fADE with constant drift and
diffusion coefficient was fitted to the same data. With constant
coefficient functions, the resulting solution process $Y_t$ in the SDE
\eqref{eq:mainSD} itself is a stable process. %The observed
%concentrations simply was used to determine the parameters of this
%stable process.
The parameters of the fitted stable process were estimated using the
least-square method.

Here we are considering an fADE that may have variable coefficients.
For the fits presented in this section, the diffusion coefficient is
still assumed to be location invariant, but we use a piecewise linear
drift function of the form \eqref{eq:Xu24}. The choice of $x_l, \
x_m$ and $x_r$ are subjective but should comply with the boundary
conditions. $x_l=0$ (starting location), $x_r=300$ (beyond the observed
location range) and $x_m=9.375$ (a rough mid-location of the observed
range) were set for the current simulation. We used the fitted
parameter values from the constant coefficient model \cite
{Chakraborty_Meerschaert} as initial values for the simulation; these
parameter values are included in Appendix~\ref{appendix}. Then we used the data to estimate the linear drift and to
obtain $p(x,t)$ that solves \eqref{eq:Xu1} by the methodology
presented in Section~\ref{s_method}. $c(x,t)=Kp(x,t)$ gives the fitted
concentration that solves fADE associated with \eqref{eq:Xu1}.

The concentrations fitted by the proposed method for the MADE site data
at day 224 and day 328 are included here.%\vadjust{\vfill{\eject}}

\medskip

\textbf{Fitted parameters for the MADE site data at day 224:}\hfill\break
$a_0=0.110$, $a_1=0.00032$, $a_2=0.0003$, $a_3=0.00019$, \hfill\break
$b= 0.1859783$, $\lambda=0.80$, $\gamma=0.9999$, $K=56778.24$.
%
%t1 #&#
\begin{table}[h]
\caption{Observed and fitted tracer concentrations for MADE data at day 224}
\label{Table_224Fit}
\begin{tabular}{|r|r|r|}
\hline
\multicolumn{3}{|c|}{Day 224}\\
\hline
\multicolumn{3}{|c|}{$x$ $=$ \mbox{location}, $C$ $=$ \mbox{tracer concentration}}\\
\hline
$x$ & Observed $C$& Fitted $C$\\
\hline
2.1000 &3378.0000 &2871.7248\\
2.9000 &1457.0000 &3187.0093\\
3.6000 &6494.0000 &3529.6230\\
6.0000 &1335.0000 &4027.3618\\
6.6000 &7705.0000 &3321.0305\\
6.8000 &2206.0000 &3086.6518\\
6.9000 &1291.0000 &2969.4625\\
7.2000 &4515.0000 &2734.3773\\
8.0000 &3598.0000 &2361.1656\\
8.7000 &2447.0000 &2083.3103\\
9.7000 &2831.0000 &1809.0524\\
10.8000 &2208.0000 &1586.8821\\
12.7000 &849.0000 &1273.2312\\
13.5000 &2213.0000 &1169.6987\\
16.0000 &1485.0000 &935.2569\\
27.5000 &443.0000 &442.0004\\
31.8000 &165.0000 &356.0560\\
40.3000 &291.0000 &245.6054\\
48.8000 &237.0000 &178.6647\\
57.5000 &76.0000 &134.2397\\
66.3000 &54.0000 &103.7457\\
74.8000 &137.0000 &82.8644\\
83.4000 &37.0000 &67.3318\\
100.7000 &28.0000 &46.5417\\
117.0000 &51.0000 &34.4212\\
166.6000 &18.0000 &16.5498\\
183.5000 &28.0000 &13.9355\\
218.5000 &9.0000 &11.4791\\
268.7000 &6.0000 &10.1389\\
\hline
\end{tabular}
\end{table}\vadjust{\eject}

\textbf{Fitted parameters for day 328 MADE site data:}\hfill\break
$a_0=0.105$, $a_1=0.00030$, $a_2=0.0005$, $a_3=0.00018$,\hfill\break
$b=0.2233695$, $\lambda=0.79$, $\gamma=0.9999$, $K=37195.05$.
%\begin{center}
%
%t2 #&#
\begin{table}[h]
\caption{Observed and fitted tracer concentrations for MADE data at day 328}
\label{Table_328Fit}
\begin{tabular}{|r|r|r|}
\hline
\multicolumn{3}{|c|}{Day 328}\\
\hline
\multicolumn{3}{|c|}{$x$ $=$ \mbox{location}, $C$ $=$ \mbox{tracer concentration}}\\
\hline
$x$ & Observed $C$& Fitted $C$\\
\hline
2.1000 &1762.0000 &2001.9586\\
3.6000 &4157.0000 &3264.5572\\
5.5000 &2018.0000 &1982.1889\\
6.0000 &1295.0000 &1653.7130\\
6.6000 &2742.0000 &1436.3243\\
6.8000 &798.0000 &1405.9525\\
6.9000 &661.0000 &1390.7665\\
7.2000 &4721.0000 &1323.8183\\
8.0000 &2877.0000 &1132.2689\\
8.7000 &946.0000 &1050.1787\\
9.7000 &3140.0000 &920.6741\\
10.8000 &2075.0000 &819.8101\\
12.7000 &1636.0000 &689.5556\\
13.5000 &1569.0000 &646.4434\\
16.0000 &869.0000 &536.8714\\
16.3000 &1825.0000 &525.8232\\
17.0000 &286.0000 &501.7569\\
18.0000 &597.0000 &470.2962\\
27.5000 &127.0000 &281.0554\\
31.8000 &49.0000 &231.1571\\
40.3000 &108.0000 &163.9576\\
48.8000 &75.0000 &121.3509\\
57.5000 &95.0000 &92.2235\\
66.3000 &38.0000 &71.8197\\
74.8000 &126.0000 &57.6513\\
83.4000 &54.0000 &47.0086\\
100.7000 &19.0000 &32.6278\\
117.0000 &46.0000 &24.1752\\
166.6000 &25.0000 &11.6326\\
183.5000 &22.0000 &9.4787\\
218.5000 &11.0000 &7.0942\\
268.7000 &8.0000 &6.0752\\
\hline
\end{tabular}
\end{table}

%f1 #&#
\begin{figure}[!ht]
\includegraphics{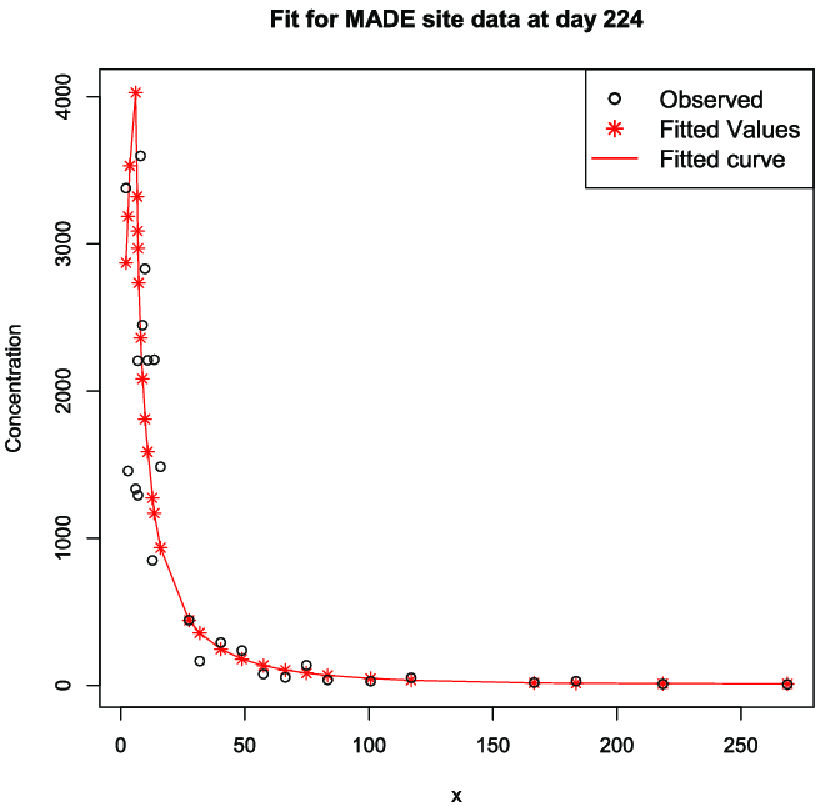}
\caption{The fitted and observed concentration at day 224 for MADE
site data}
\label{fig: 224_fit}
\end{figure}

%f2 #&#
\begin{figure}[!ht]
\includegraphics{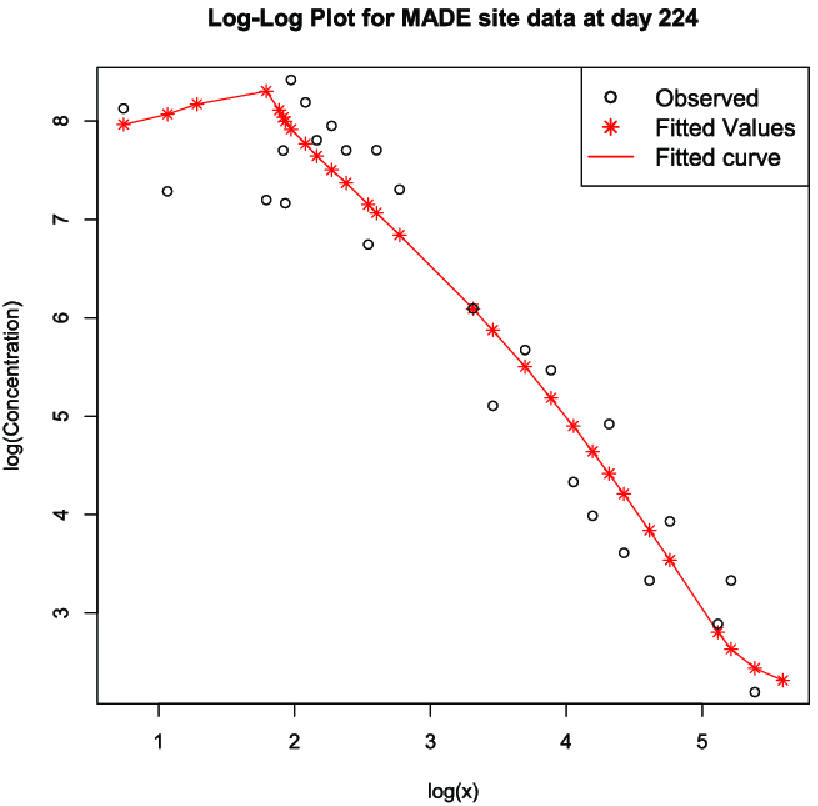}
\caption{The fitted and observed concentration with log scale at day
224 for MADE site data}
\label{fig: 224_log}
\end{figure}

%f3 #&#
\begin{figure}[!ht]
\includegraphics{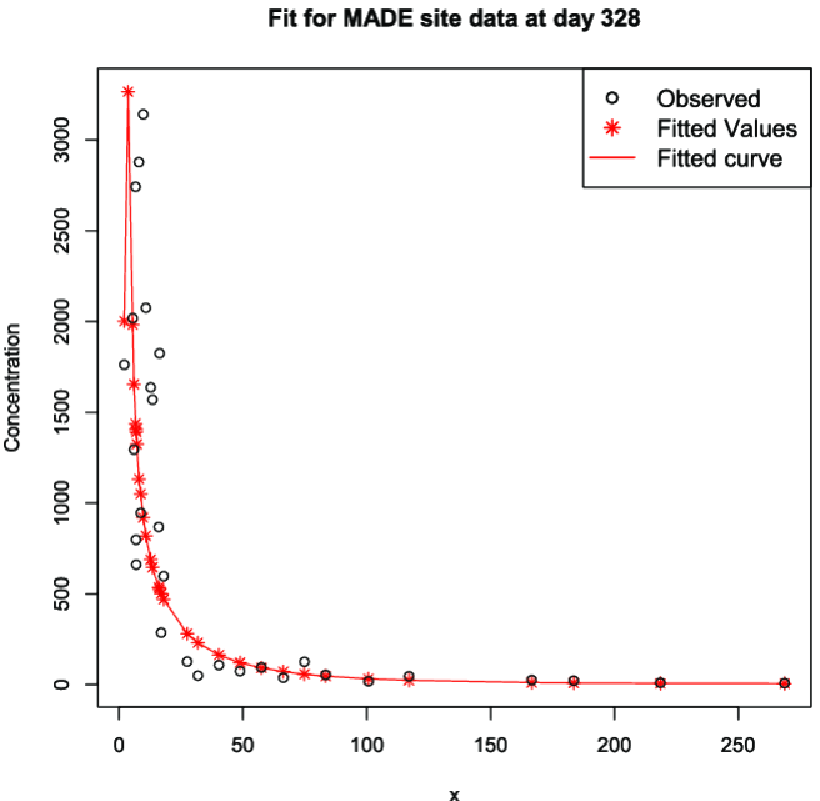}
\caption{The fitted and observed concentration at day 328 for MADE
site data}
\label{fig: 328_fit}
\end{figure}

%f4 #&#
\begin{figure}[!ht]
\includegraphics{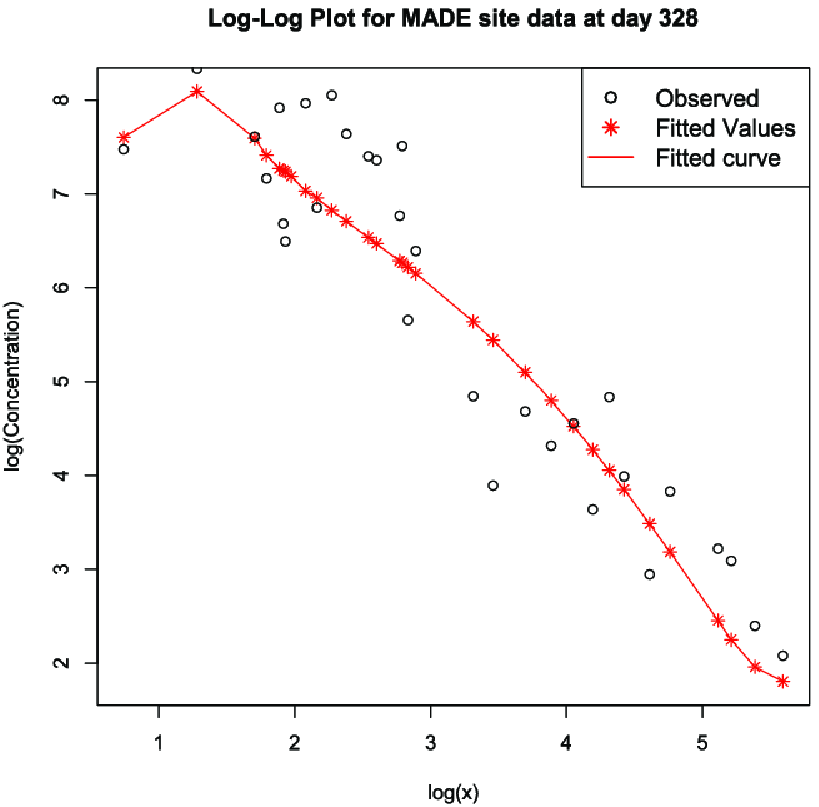}
\caption{Fitted and observed concentration with log scale at day 328
for MADE site data}
\label{fig: 328_log}
\end{figure}

\begin{figure}[!ht]
\includegraphics{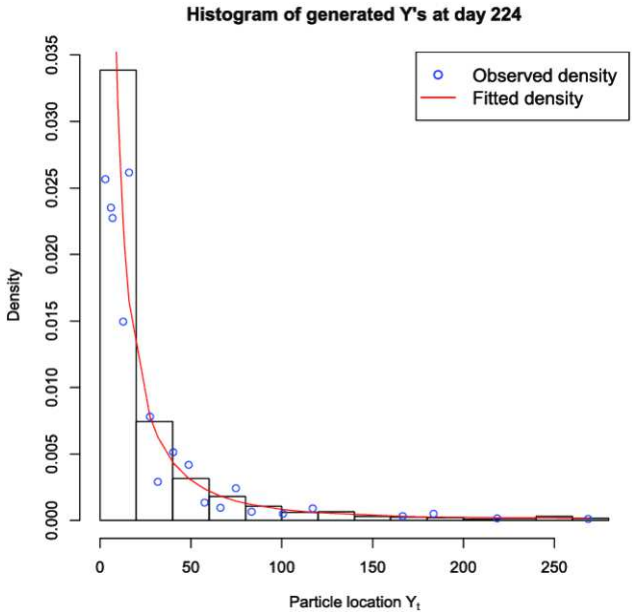}
\caption{Histogram of 1000 $Y_t$ values generated by the proposed
method at $t=224$. The red line shows the fitted pdf from \eqref
{eq:Xu1} and blue dots show observed density = observed
concentration\,/$K$ (data estimates $K=56778.24$)}
\label{fig: hist224}
\end{figure}

\noindent\textbf{Comments}:
\begin{enumerate}
\item The fit above uses piecewise linear drift velocity and location
invariant diffusion coefficient as opposed to constant drift but same
diffusion coefficient used in \cite{Chakraborty_Meerschaert}. The
fitted curves in Figures~\ref{fig: 224_log} and~\ref{fig: 328_log}
show a better tail fit here compared to similar plots in \cite
{Chakraborty_Meerschaert}. The fractional order $\alpha$ for the fADE
with constant coefficient was fitted as $1.09$ and $1.05$ for day 224
and day 328 respectively in \cite{Chakraborty_Meerschaert}. Here the
fitted fractional order $\alpha(=2-\lambda)$ are $1.2$ and $1.19$ for
day 224 and day 328 respectively. Since the new fit improves on the fit
in \cite{Chakraborty_Meerschaert}, this indicates that a better drift
velocity fit may lead to an adjustment of the fractional order to
improve the tail fit in the fADE. Ideally, the drift velocity function
might be mechanically estimated from a properly designed field
experiment while the diffusion coefficient and the fractional order can
be estimated as described in this paper for an even better fADE fit.
\item Using simulation step 2, iid sample of size 1000 for $Y_t$ was
generated at $t=\mbox{day }224$ and $t=\mbox{day }328$. Histograms
for the generated samples in Figures~\ref{fig: hist224} and~\ref{fig:
hist328} are fairly consistent with the observed concentration (or
observed density) data. These figures indicate that although we had
adopted a numerical approach to generate $Y_t$'s, the simulation method
is efficient enough to follow the underlying probability distribution.
\end{enumerate}

%f5 #&#

%f6 #&#
\begin{figure}[!ht]
\includegraphics{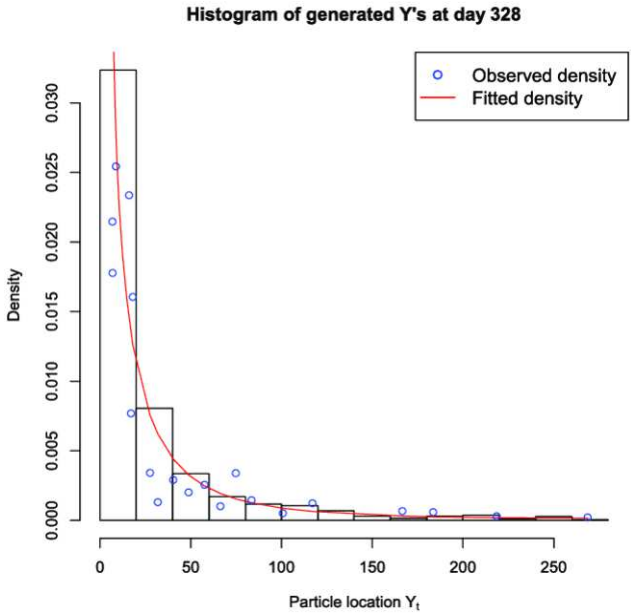}
\caption{Histogram of 1000 $Y_t$ values generated by the proposed
method at $t=328$. The red line shows the fitted pdf from \eqref
{eq:Xu1} and blue dots show observed density = observed
concentration\,/$K$ (data estimates $K=37195.05$)}
\label{fig: hist328}
\end{figure}

%s4 #&#
\section{Discussion}\label{sec4}
The trapezoidal rule provides a very simple and quick numerical
approximation of the CDF of $Y_t$ in \eqref{eq:mainSD} and for its
inverse in the proposed simulation step 2. A more sophisticated
numerical method like Simpson's rule or other higher order
Newton--Cotes formula or quadrature rules \cite{Scarborough} can be
used to obtain a more accurate $\hat{F}_{Y_t}$, but the procedure will
be much more computationally taxing.

The method described here can be used to simulate a large number of
observations (approximated) from $Y_t$. With the usual particle
tracking insights, these simulated values can also be used to build
empirical confidence intervals for the fitted concentration values.

Let us consider $P[Y_t\in(a,b)]=F_{Y_t}(b)-F_{Y_t}(a)$. This is the
probability that a randomly chosen tracer particle will be in a given
interval $(a,b)$ at a time point $t$. For example, in groundwater
pollution modeling this expression can be an important one that
estimates the chance that the pollution will reach a certain area after
a certain time. Using i.i.d. random simulations $Y^{(1)}_t,
Y^{(2)}_t,\ldots, Y^{(N)}_t $ of $Y_t$, an empirical estimate for this
probability can be given by $\hat{P}=\frac{1}{N}\sum^N_{i=1}I[Y^{(i)}_t\in(a,b)]$. Then by the central limit theorem $\hat
{P}$ asymptotically follows the normal distribution with mean $P[Y_t\in
(a,b)]$ and variance $\dfrac{P[Y_t\in(a,b)](1-P[Y_t\in(a,b)])}{N}$.
Hence for large $N$'s, an asymptotic $(1-\alpha)100\%$ confidence
interval for $P[Y_t\in(a,b)]$ can be given by
%
%e37 #&#
\begin{equation}
\label{eq:CI1} \hat{P}\mp z_{\alpha/2}\sqrt{\frac{\hat{P}(1-\hat{P})}{N}}
\end{equation}
Now noting that density $p(y,t)$ can be approximated by $\frac
{1}{b-a}P[Y_t\in(a,b)]$ with a small interval $(a,b)$, the simulated
$Y_t$ values can be used to calculate the confidence intervals for the
fitted densities. See \cite{Chakraborty_Meerschaert} for detailed
asymptotic confidence interval construction steps and related results
associated with the empirical density function that is calculated using
the generated $Y_t$'s.

Assuming that the observed concentrations are scaled versions of
density $p(y,t)$ of $Y_t$, one can repeat simulation step 2 described
in Section~\ref{subs:simulation} without solving the forward equation
to generate $Y_t$'s. However, the grid used for the simulation will be
limited only to the fixed observed data location points. Since $Y_t$'s
are simulated through numerical approximation the performance of this
simulation depends on the grid spacing used. Thus the fitted fADE
%%equation
in simulation step 1 not only gives us the better understanding of the
underlying process but also is essential for generating a good
approximate process that resembles $Y_t$ by enabling the use of a finer
grid for the numerical approximation. Further, modeling the underlying
fADE process facilitates better prediction than mere extrapolation from
an observed sample set.

On the other hand, an fADE can be fitted to the observed concentration
by simulation step 1 only, without describing the underlying stochastic
process. But the stochastic diffusion description is useful for
understanding the mesoscopic flow and particle tracking methods.
Further, generating the solution to the SDE associated with the fADE
can be used to construct confidence intervals for the fADE fits that
account for the error of estimation by a sample set.

In conclusion, the fADE and the associated $\alpha$-stable SDE are
essential tools to model heavy-tailed diffusion. However, the model
fitting part may get complicated if we consider any scenario other than
constant drift and diffusion coefficients. A numerical scheme for
solving the fADE and the related SDE with linear drift is presented
here. The general ideas used in the simulation steps can be replicated
for a more complicated form of the drift. We plan to explore such
models in our future work along with the convergence and the efficiency
of the numerical approximations applied here.

%\bibliographystyle{abbrv}
%\bibliography{SDE_bib}

%sA #&#
\begin{appendix}
\appendix
%sB #&#
\section{Appendix}\label{appendix}

In \cite{Chakraborty_Meerschaert} the tracer trajectory $Y_t$ of the
MADE site data was assumed to follow an SDE of the form \eqref
{eq:mainSD} with constant parameters and hence itself was a stable
process with index of stability $\alpha$, skewness parameter $\beta$,
location parameter $v$ (same as the assumed constant drift velocity in
\eqref{eq:fADE}) and scale parameter $\sigma$, where $\sigma^\alpha
=Dt|\cos(\pi\alpha/2)|$ with constant diffusion coefficient function
$D$ in \eqref{eq:fADE}.

The fitted values in \cite{Chakraborty_Meerschaert} are:

\smallskip

\textbf{Day 224 fitted values:}\hfill\break
$\sigma=5.137167$, $\mu=43.915430$;\hfill\break
$\beta=0.99$;\hfill\break
$\alpha=1.0915$;\hfill\break
$v=0.196051 m/\textit{day}$\hfill\break
$D= 0.1859783 m^\alpha/\textit{day}$\hfill\break
$K= 56778.24$.

\smallskip

\textbf{Day 328 fitted values:}\hfill\break
$\sigma=5.380654 $; $\mu=74.677975$;\hfill\break
$\beta=0.99 $;\hfill\break
$\alpha=1.050998 $;\hfill\break
$v= 0.2276768 m/\textit{day}$\hfill\break
$D= 0.2233695 m^\alpha/\textit{day}$\hfill\break
$K= 37195.05$%\\

\smallskip

These were used as initial values for the iterations in Section~\ref{s:illustration}.
\end{appendix}

%\begin{acknowledgement}%[title={Acknowledgments}]
%\end{acknowledgement}

%\begin{funding}
%\gsponsor[id=,sponsor-id=]{}
%\gnumber[refid=]{}
%\end{funding}

%
\end{document}